%% file: affine-rf-main.tex
\documentclass[11pt]{amsart}

\usepackage{graphicx}
\usepackage{amscd}
\usepackage{amsmath}
\usepackage{amsfonts}
\usepackage{amssymb}
\usepackage{amsxtra}
\usepackage{xspace}
\usepackage{array}

\theoremstyle{plain}
\newtheorem{theorem}{Theorem}

\newtheorem{corollary}{Corollary}
\newtheorem{proposition}{Proposition}

\theoremstyle{definition}

\numberwithin{equation}{section}

\newcommand{\be}{\begin{enumerate}}
\newcommand{\ee}{\end{enumerate}}
\newcommand{\beq}{\begin{equation}}
\newcommand{\eeq}{\end{equation}}
\newcommand{\bprop}{\begin{proposition}}
\newcommand{\eprop}{\end{proposition}}

\newcommand{\reals}{\mathbb{R}}

\newcommand{\rat}{\mathbb{Q}}
\newcommand{\integers}{\mathbb{Z}}

\newcommand{\pfbegin}{\noindent {\em Proof:} }

\DeclareMathOperator{\des}{des} 
\DeclareMathOperator{\len}{\ell}
\DeclareMathOperator{\refl}{refl}

\newcommand{\refw}[1][\afw]{\refl(#1)}
\newcommand{\highroot}{\widetilde{\alpha}} 

\newcommand{\p}{^{\prime}}

\newcommand{\ps}{Poincar\'{e} series\xspace}
\newcommand{\rfs}{rational functions\xspace}

\newcommand{\afw}{\widetilde{W}} 
\newcommand{\afs}{\widetilde{S}} 
\newcommand{\afv}{\widetilde{V}} 
\newcommand{\afrs}[1][]{\widetilde{\Phi}^{#1}}  

\newcommand{\wta}[1][A]{\afw^{#1}}  

\newcommand{\wto}{\afw^{0}}  

\newcommand{\wtxu}{\afw_{x,u}}  
\newcommand{\wtaxu}{\afw^A_{x,u}}

\newcommand{\fw}{W}  
\newcommand{\frs}[1][]{\Phi^{#1}}  

\newcommand{\fcrs}[1][]{\Check{\frs}^{#1}}  

\newcommand{\qc}[1][]{\Check{Q}^{#1}}  

\newcommand{\xwp}{X_{W\p}}  

\newcommand{\alc}{A_f}

\begin{document}

\title[]{Poincar\'{e} series of subsets of affine Weyl groups}
\author{Sankaran Viswanath}
\address{Department of Mathematics\\
University of California\\
Davis, CA 95616, USA}
\email{svis@math.ucdavis.edu}

\begin{abstract}
In this note, we identify a natural class of subsets of affine Weyl
groups whose Poincar\'{e} series are rational functions. This class
includes the sets of minimal coset representatives of reflection
subgroups. As an application, we construct a generalization of the
classical length-descent generating function, and prove its rationality.
\end{abstract}

\maketitle

\input{affine-rf.tex}

\bibliographystyle{amsplain}
\bibliography{biblio-affine-rf,biblio-cox-exp}
\end{document}

%% file: affine-rf.tex
\section{Introduction}
The aim of this note is to prove the rationality of certain length
generating functions in affine Weyl groups.
Let $\afw$ be an affine Weyl group with Coxeter generators $\afs =\{s_i\}_{i=0}^n$ and length function $\len(\cdot)$. 
 For a subset $X$ of $\afw$, its
Poincar\'{e} series (length generating function) is $X(q) := \sum_{w
  \in X} q^{\len(w)}$. For many ``natural'' subsets $X$, $X(q)$ turns out
to be a rational function; notable examples include the entire group $\afw$, its parabolic
subgroups $\afw_I$ ($I \subset \afs$), and the sets $\afw^I$ of minimal length left coset representatives
for $\afw_I$ in $\afw$ (in fact all this holds for arbitrary Coxeter groups $\afw$).

One of our objectives is to study an interesting class of subgroups of $\afw$ -
{\em reflection subgroups}. These are subgroups  $W\p \subset \afw$
generated by reflections (conjugates of elements of $\afs$); they
exhibit many of the same properties as parabolic subgroups
\cite{deodhar2}, \cite{dyer1} - for e.g (i) $W\p$ is a Coxeter group in its own right, (ii) $\exists$ unique minimal length left coset representatives for $W\p$ in $\afw$, etc. A complete classification of 
the reflection subgroups of affine Weyl groups (in terms of those of the underlying finite Weyl group) was given by Dyer in \cite{dyer1}.

However, reflection subgroups are ill-behaved wrt the length
function. The length function of $W\p$ does not in general agree with the restriction of $\len(\cdot)$
to $W\p$. This makes it difficult to study the
Poincar\'{e} series $W\p(q)$; for instance, it does not seem to be
known if $W\p(q)$ is a rational function for all reflection subgroups
$W\p$. Our focus however, will be on the set  $\xwp$ of
 minimal length left coset representatives for $W\p$; a complication here
is that $\xwp(q) \neq \afw(q)/W\p(q)$ any more (this holds if $W\p$ is
parabolic). This means that even in cases where $W\p(q)$ is known to
be a rational function (e.g $\# W\p < \infty$), we still cannot
conclude that $\xwp(q)$ is rational. Our first goal is to show that
$\xwp(q)$ is indeed a rational function for all reflection subgroups
$W\p$ of $\afw$. We thus get another natural
class of subsets of $\afw$ with rational \ps. We remark that if $\afw$
is allowed to be an arbitrary (non finite, non affine) Coxeter group,
the rationality question for $\xwp(q)$ seems harder to decide; the
article \cite{vis} is concerned with reflection subgroups of such
Coxeter groups, but only deals with the {\em growth type} of the set $\xwp$.

Next, we turn to a two variable refinement of the \ps of $\afw$  - the classical
length-descent generating function $\afw(q,t):=\sum_{w \in \afw}
q^{\len(w)} t^{\,\des(w)}$ where $\des(w):=\#\{s \in \afs: \len(ws) <
\len(w)\}$. This is well known (see \cite{reiner}) to be a polynomial in
$t$ with coefficients that are \rfs in $q$, i.e $\afw(q,t) \in
\rat(q)[t]$. To generalize this, let  $\refw := \displaystyle\bigcup_{\sigma \in \afw} \sigma \afs
\sigma^{-1}$ be the set of reflections in $\afw$ and $A \subset \refw$
be a finite subset.  Define
$$ \afw(q,t,A):= \sum_{w \in \afw} q^{\len(w)} t^{\,\des_A(w)}$$
 where $\des_A(w):=\#\{r \in A: \len(wr) < \len(w)\}$. Thus
 $\afw(q,t,\afs) = \afw(q,t)$. 

\vspace{0.1cm}
\noindent
The second aim of this note is to show
 that $\afw(q,t,A) \in  \rat(q)[t]$ for all finite $A \subset \refw$.

Both this and the earlier result on reflection subgroups will be shown
to fit into a slightly more general framework. They will follow as simple consequences
of our main theorem (theorem \ref{mainthm}) which also seems to   be of
independent interest.

{\em Acknowledgements:} I'd like to thank Monica Vazirani for helpful
discussions while this work was in progress.

\section{The main theorem} \label{two}
\subsection{Preliminaries}
Let $\fw$ be the finite Weyl group corresponding to an irreducible,
crystallographic root system $\frs$. Let $\{s_i\}_{i=1}^n, \{\alpha_i\}_{i=1}^n,
 \{\Lambda_i\}_{i=1}^n$ be the simple reflections, the simple roots
 and the fundamental weights respectively.  Let $V$ be the
 $\reals$ span of the simple roots and $(,)$ be a positive
 definite, $\fw$ invariant bilinear form on $V$. Similarly, one has the
 coroot system $\fcrs \subset V$ with simple coroots $\Check{\alpha}_i$ and
 fundamental coweights $\Check{\Lambda}_i$. The root and coroot
 lattices will be denoted $Q = \integers \frs$ and $\qc =
 \integers \fcrs$. 
For $\alpha \in V$, let $t_\alpha$ be the translation map  $v
\mapsto v + \alpha$. Let $T = T(\qc) := \{t_\alpha : \alpha \in \qc\}$
be the group of translations of $V$ by elements of $\qc$. The affine Weyl
group $\afw$ can be defined as the subgroup of the group of affine transformations
of $V$ generated by $\fw$ and $T$; we have $\afw = W \ltimes T$. 

$\afw$ is a Coxeter group with generators $\afs:= \{s_i\}_{i=1}^n \cup
\{s_0\}$; here $s_0$ is the reflection about the affine hyperplane
$\{v \in V: (v,\highroot)=1\}$ with $\highroot =$ the highest root of
$\frs$. Thus, $\afw$  also has a {\em reflection representation} (or geometrical realization)
$\afv$. To construct this (see \cite[\S 4]{jrs} or \cite[Chap 6]{kac}),
 set $\afv := V \oplus \reals \delta$ and
extend $(,)$ to a positive semidefinite form on $\afv$ by letting $(\delta, \delta)=0, (\delta,v) =0
\, \forall v \in V$.  The $\afw$ action on $\afv$ extends the $\fw$
action on $V$ via the following prescription: given $\sigma \in \fw, \alpha \in
\qc$, $v \in V$,
\begin{align*}
\sigma \delta &= t_{\alpha} \delta = \delta\\
 t_\alpha(v) &= v - (v,\alpha) \,\delta 
\end{align*}
The root system of  $\afw$ is $\afrs :=\{\beta + k \delta: \beta \in
\frs, k \in \integers\} \subset \afv$; the  simple roots are $\{\alpha_i\}_{i=0}^n$,
where $\alpha_0 = \delta - \highroot$. The positive roots of $\afw$ are :
\beq \label{posroots}
\afrs[+] = \{\beta + k \delta: \beta \in \frs[+], k \geq 0\} \cup 
\{\beta + k \delta: \beta \in \frs[-], k \geq 1\}
\eeq
Given $\gamma \in \afrs[+]$, one has the map $s_\gamma \in GL(\afv)$ defined by  $s_\gamma(\mu) := \mu -
2\frac{(\mu,\gamma)}{(\gamma, \gamma)} \,\gamma$. It is well known
that the set $\{s_\gamma: \gamma \in \afrs[+]\}$ is the image of
$\refw$ in $GL(\afv)$; thus $\afrs[+]$ is in bijection with $\refw$.

\subsection{The main theorem and its corollaries}
Let $\len(.)$ be the length function on $\afw$ wrt $\afs$
(this extends the length function on $W$). Given a finite subset $A
\subset \refw$, let $\wta :=\{ \sigma \in \afw : \len(\sigma r) >
\len(\sigma) \; \forall r \in A\}$. Observe that if $A\p$
is the corresponding set $\{\gamma \in \afrs[+]: s_\gamma \in A\}$, we
have the equivalent definition $\wta = \{  \sigma \in \afw :
\sigma(A\p) \subset \afrs[+]\}$. We also note that in the 
familiar case when $A = I \subset \afs$, $\wta$ is just the set
of minimal left coset representatives for the parabolic subgroup
$\afw_I$. Our main theorem is :
\begin{theorem}\label{mainthm}
For any finite $A \subset \refw$, the \ps 
$\wta(q) = \displaystyle\sum_{w \in \wta} q^{\len(w)}$ is a rational function.
\end{theorem}
We postpone the proof to section \ref{pfsection}. We first use this
theorem to give quick proofs of the two results mentioned in the introduction.

\subsubsection{}
If $W\p \subset \afw$ is a reflection subgroup, it is a well known
theorem due (independently) to Deodhar \cite{deodhar2} and Dyer
\cite{dyer1} that $W\p$ is a Coxeter group wrt a set $S\p =
\{s_{\gamma_i}:i=1\cdots k\}$ of reflection generators. Here, the
$\gamma_i \in \afrs[+]$ and satisfy $(\gamma_i, \gamma_j) \leq 0 \;
\forall i \neq j$ \cite[Theorem 4.4]{dyer1} (this theorem holds even
when $\afw$ is an arbitrary Coxeter group in which case $S\p$ need not
be finite; for affine $\afw$ however, it is easy to show that $\# S\p
< \infty$). In \cite{dyer1}, Dyer also showed that there are unique minimal length
elements in the left cosets of $W\p$. Let $\xwp$ denote this set of
minimal coset representatives; then  $\sigma \in \xwp \Leftrightarrow
\len(\sigma s_{\gamma_i}) > \len(\sigma) \; \forall i=1\cdots k$. Thus
$\xwp = \wta[S\p]$. As a consequence of theorem \ref{mainthm},we have
\begin{corollary}
Let $W\p$ be any reflection subgroup of $\afw$. Then $\xwp(q)$ is a
rational function.
\end{corollary}

\subsubsection{}
We recall from the introduction that the generating function
$\afw(q,t) = \sum_{w \in \afw} q^{\len(w)} t^{\des(w)} \in
\rat(q)[t]$. We refer to Reiner's article \cite[Theorem 1]{reiner} for the proof
of this ``folklore'' result. The proof essentially consists in showing the
following identity (in our notation):
\beq \label{reiner1}
\afw(q,t) = \sum_{I \subset \afs} t^{|I|} \,(1-t)^{|\afs \,\backslash I|}\;
\wta[\afs \,\backslash I] (q)
\eeq
 Here $|\cdot|$ denotes set cardinality. To complete Reiner's
 argument, one observes that since $\afs \,\backslash I \subset \afs$,
 we have $\wta[\afs \,\backslash I]  =  \dfrac{\afw(q)}{\afw_{\afs \,\backslash I}(q)}$,
 which is a rational function. 

It is now elementary to modify the above argument  for
the case where $\afs$ is replaced by $A$. The analogue to equation
\ref{reiner1} is now :
\beq
\afw(q,t,A) = \sum_{B \subset A} t^{|B|} \, (1-t)^{|A    \,\backslash B|}
\; \wta[A \,\backslash B] (q)
\eeq
where $\afw(q,t,A):= \sum_{w \in \afw} q^{\len(w)} \, t^{\,\des_A(w)}$.
 Invoking  theorem \ref{mainthm}, we conclude $\wta[A
  \,\backslash B] (q) \in \rat(q)$ and hence the following:
\begin{corollary}
 $\afw(q,t,A) \in \rat(q)[t]$ for all finite subsets $A \subset \refw$.
\end{corollary}

\section{Proof of main theorem} \label{pfsection}
\subsection{}
Before embarking on the proof of our main theorem, we collect together
some well known facts concerning $\afw$ (good references are
\cite{humphreys}, \cite{bourbaki}). We freely use the notation of section
\ref{two}. Let 
$$C = \{ v \in V: (v, \alpha_i) \geq 0 \, \forall i=1\cdots n\}; \;\; 
\alc = \{ v \in C: (v, \highroot) \leq 1\}$$
be the closures of the fundamental chamber and fundamental alcove respectively. The
finite Weyl group $\fw$ is a parabolic subgroup of $\afw$; let $\wto$ be the set of minimal
length right coset representatives for $\fw$ in $\afw$.

\vspace{0.2cm}
\noindent
{\bf  Fact 1:} $\forall w \in \afw$, $\exists ! u \in \fw$ s.t $uw \in
\wto$; this $u$ is the unique element of $\fw$ s.t $uw(\alc) \subset
C$. 

\vspace{0.2cm}
Let $\rho \in V$ be the Weyl vector of $\frs$; it is determined by the conditions
$(\rho, \Check{\alpha}_i) = 1 \; \forall i=1\cdots n$.

\vspace{0.1cm}
\noindent
{\bf Fact 2:} $\len(t_\alpha) = \len(t_{\sigma \alpha}) \; \forall
\sigma \in \fw, \alpha \in \qc$; further $\len(t_\alpha) = (\alpha,
2\rho)$ if $\alpha$ is a dominant element of $\qc$.

\vspace{0.2cm}
For $u \in \fw$, define $T_u :=\{t_\alpha \in T: ut_\alpha \in
\wto\}$. By fact 1, this means $u t_\alpha (A_f) \subset C$.
Let the  vertices of the  simplex $\alc$ be $\{0, \theta_1, \theta_2, \cdots, \theta_n\}$. 
We note that the $\theta_j \in C$, but in general they are not elements of
the coweight lattice; we only have $(\theta_j, \alpha_i) \in \rat \;
\forall i,j$. It is clear that $ut_\alpha(A_f) \subset C \Leftrightarrow u\alpha \in C
\text{ and } u(\alpha + \theta_j) \in C \; \forall j=1\cdots n$. Let
$m_i(u)$ be the smallest integer such that $m_i(u) \geq 0$ and $m_i(u) \geq
-(u \theta_j, \alpha_i) \; \forall j=1\cdots n$. Since $u\alpha \in
\qc$, the above discussion implies 

\vspace{0.1cm}
\noindent
{\bf Fact 3:} The condition $t_\alpha \in T_u$ is equivalent to the system of inequalities : 
\beq \label{ineq1}
(u\alpha, \alpha_i) \geq m_i(u) \;\; \forall i=1\cdots n
\eeq

\subsection{Proof of theorem \ref{mainthm}}
We refer back to the statement of  theorem \ref{mainthm}. We will prefer to
work with the set $A\p = \{\gamma \in \afrs[+]: s_\gamma \in A\}$
rather than with $A$ itself. Thus $\wta = \{ \sigma \in \afw:
\sigma(A\p) \subset \afrs[+]\}$. First, we can assume wlog
that for each $\beta \in \frs$, $A\p$ contains at most one element of
the form $\beta + k \delta$. If not, suppose $\beta + k_1 \delta,
\, 
\beta + k_2 \delta \in A\p$ with $k_1 < k_2$. For $\sigma \in \afw$,
$\sigma(\beta + k_1 \delta) = \sigma \beta + k_1 \delta \in \afrs[+]
\; \Rightarrow \sigma \beta + k_2 \delta \in \afrs[+]$ too. Thus we
can delete $\beta + k_2 \delta$ from $A\p$ without changing $\wta$.

For each $\beta \in \frs$, let $k_\beta$ be the unique integer (if it
exists) such that $\beta + k_\beta \,\delta \in A\p$. If $A\p$ contains
no element of the form $\beta + k \,\delta$, we set $k_\beta
:=\infty$. Let $F:=\{\beta \in \frs: k_\beta < \infty\}$.

Next, we'll analyze what it means for $\sigma$ to be an element of
$\wta$. We write $\sigma = xt_\alpha$ with $x \in \fw$, $\alpha \in
\qc$. For each $\beta \in F$, we require
$$\sigma (\beta + k_\beta \,\delta) = xt_\alpha(\beta + k_\beta \,\delta) =
x \beta + (k_\beta - (\alpha, \beta))\,\delta \in \afrs[+]$$
By equation \ref{posroots}, this implies that $\alpha$ satisfies the
following inequalities 
\beq \label{ineq2}
(\alpha, \beta) \leq  
\begin{cases}
k_\beta & \text{ if } x\beta \in \frs[+] \\
k_\beta  -1 & \text{ if } x\beta \in \frs[-] 
\end{cases}
\eeq
Note that by our convention of setting $k_\beta = \infty$ for $\beta
\not\in F$, we can state this as: $xt_\alpha \in \wta \Leftrightarrow$
 the inequalities  \eqref{ineq2}  hold for {\em all} $\beta \in \frs$
 (not just for $\beta \in F$).

Now, for fixed $x,u \in \fw$, define $\wtxu:=\{xt_\alpha: t_\alpha \in
T_u\}$ and $\wtaxu := \wta \cap \wtxu$. Then $\wta = \displaystyle\bigsqcup_{x,u \in \fw} \wtaxu$. Given
$\sigma \in \wtxu$, we have $\sigma = xt_\alpha =
(xu^{-1})(ut_\alpha)$; since $ut_\alpha \in \wto$, $$\len(\sigma) =
\len(xu^{-1}) + \len(ut_\alpha) = \len(xu^{-1}) + \len(t_\alpha) -
\len(u)$$
Thus  
\beq \label{maineq}
 \sum_{\sigma \in \wtaxu} q^{\len(\sigma)} = q^{\len(xu^{-1}) -
  \len(u)} \sum_{xt_\alpha \in \wtaxu} q^{\len(t_\alpha)}
\eeq
Let $ f_{x,u}(q):= \displaystyle\sum_{xt_\alpha \in \wtaxu} q^{\len(t_\alpha)}$; by
fact 2, we have $ f_{x,u}(q) = \displaystyle\sum_{xt_\alpha \in \wtaxu}
q^{\len(t_{u\alpha})}$.

\vspace{0.2cm}
\noindent
{\em Claim:} $f_{x,u}(q)$ is a rational function.

\pfbegin Observe that $xt_\alpha \in \wtxu$ iff $\alpha$ satisfies the inequalities
\eqref{ineq1}  and $xt_\alpha \in \wta$ iff $\alpha$ satisfies the
inequalities \eqref{ineq2}. Now, since $(\alpha,\beta) =
(u\alpha,u\beta) \, \forall \beta \in \frs$, inequalities \eqref{ineq2}
can be rewritten (with $\gamma = u\beta$) as : 
\beq \label{ineq3}
\forall \gamma \in \frs, \;\;\;(u\alpha, \gamma) \leq 
\begin{cases}
k_{u^{-1} \gamma} & \text{ if } xu^{-1}\gamma \in \frs[+] \\
k_{u^{-1}\gamma}   -1 & \text{ if } xu^{-1}\gamma  \in \frs[-] 
\end{cases}
\eeq
So, $xt_\alpha \in \wtaxu$ iff $u\alpha$ satisfies the  systems of inequalities
\eqref{ineq1} and \eqref{ineq3}. Since $u\alpha \in C$, fact 2 also gives
$\len(t_{u\alpha}) = (u\alpha, 2\rho)$. Thus 
$$f_{x,u}(q) = \sum_{xt_\alpha \in \wtaxu} q^{(u\alpha, 2 \rho)} =
\sum \, q^{(\delta, 2\rho)}$$ 
where the last sum on the right runs over all $\delta
\in C \cap \qc \subset \qc[+]$ satisfying the systems of inequalities
\eqref{ineq1} and \eqref{ineq3} (with $\delta$ in place of $u\alpha$).

Now, since $\delta \in \qc[+]$ is a non-negative integer linear
combination of simple coroots, the set of allowed $\delta$ in the
above summation can be
thought of as the solution set in non-negative integers of a system of
inequalities with integer coefficients. By the classical theory of
such systems (see for e.g \cite[\S 4.6]{stanley1}, \cite{stanley2}), the generating series $f_{x,u}(q)= \sum
q^{(\delta, 2\rho)}$ is a rational function.

Finally, since $\wta(q) = \displaystyle\sum_{x,u \in \fw} q^{\len(xu^{-1}) -
  \len(u)} f_{x,u}(q)$, it is clear that $\wta(q)$ is a rational function. This
  completes the proof of our main theorem. \qed

\subsection{}
As a by-product of our  method of proof above, we obtain the following
fact concerning the rationality of the \ps
$T(q) = \sum_{t_\alpha \in T} q^{\len(t_\alpha)}$.
Since $T = \bigsqcup_{u \in W} T_u$, we have
\begin{align*}
T(q) &=  \sum_{u \in W}\sum_{t_\alpha \in T_u}  q^{\len(t_\alpha)} =  \sum_{u \in W}
\sum_{t_\alpha \in T_u}  q^{\len(t_{u\alpha})}\\
 &=  \sum_{u \in W} \sum_{t_\alpha \in T_u}  q^{(u\alpha,2\rho)}
\end{align*}
 The set $\{u \alpha: t_\alpha \in T_u\}$ is precisely the set of
 elements in $C \cap \qc$ satisfying the inequalities
\eqref{ineq1}; again by the general theory quoted above,
we conclude that $\sum_{t_\alpha \in T_u} \,q^{(u\alpha,2\rho)}$ is a
 rational function. This proves 
\begin{corollary}
 $T(q) = \sum_{t_\alpha \in T} q^{\len(t_\alpha)}$ is a rational
 function.
\end{corollary}

As our final remark, we compare the result of the above corollary 
with a related fact about $T$ that can be derived from
general considerations concerning finitely generated abelian groups.
If K is any finite set of generators of the (free) abelian group $T$,
we have the length function $\len_K(t_\alpha)$, defined to be the length of the
smallest word in $K \cup K^{-1}$ that represents $t_\alpha$. It is
well known (see for e.g \cite{benson}) that the generating series $\sum_{t_\alpha \in T} q^{\len_K(t_\alpha)}$ is a
  rational function. In our situation above however, the length function $\len$ on $T$ is wrt the Coxeter
  generators $\afs$ of the ambient group $\afw$ (note that none of
  these generators is in $T$).